\documentclass[12pt]{article}
\usepackage{amssymb,amsmath, amsthm}
\newtheorem{thm}{Theorem}
\newtheorem*{theo}{Theorem}

\newtheorem{lemma}[thm]{Lemma}
\newenvironment{defin}{\medskip\noindent{\sc
Definition}.}{\goodbreak\medskip}
\newenvironment{nota}{\medskip\noindent{\sc
Notations}.}{\goodbreak\medskip}
\newenvironment{remk}{\noindent{\sc
Remark}.}{\goodbreak\vskip10pt}
\newtheorem{prop}[thm]{Proposition}

\def\demo{\medskip\goodbreak\noindent
     \hbox{\sc Proof \kern .3em}\ignorespaces}%
  \def \qedbox{$\square$}%
  \def \qed{\hglue1mm\hfill{\ifmmode\qedbox
     \else\unskip\ \hglue0mm\hfill\qedbox\medskip
      \goodbreak\fi}}%
\def\enddemo{\qed\goodbreak\vskip10pt}%

\newcommand{\stable}{\mathbb {S}}
\newcommand{\T}{\mathbb {T}}
\newcommand{\U}{\mathbb {U}}
\newcommand{\A}{\mathbb {A}}
\newcommand{\esse}{\mathbb {S}}
\newcommand{\R}{\mathbb {R}}

\newcommand{\Z}{\mathbb {Z}}
\newcommand{\N}{\mathbb {N}}

\newcommand{\Hc}{\mathcal {H}}
\newcommand{\Cc}{\mathcal {C}}

\newcommand{\Gc}{\mathcal {G}}
\newcommand{\Lc}{\mathcal {L}}

\newcommand{\Ac}{\mathcal {A}}
\newcommand{\Tc}{\mathcal {T}}

\newcommand{\Sc}{\mathcal{S}}

\begin{document}
\title{{ Lower and upper bounds  for the Lyapunov exponents of twisting dynamics: a relationship between the exponents and the angle of the Oseledet's splitting }}
\author{M.-C. ARNAUD
\thanks{ANR Project BLANC07-3\_187245, Hamilton-Jacobi and Weak KAM Theory}
\thanks{ANR DynNonHyp}
\thanks{Universit\'e d'Avignon et des Pays de Vaucluse, Laboratoire d'Analyse non lin\' eaire et G\' eom\' etrie (EA 2151),  F-84 018Avignon,
France. e-mail: Marie-Claude.Arnaud@univ-avignon.fr}
}
\maketitle
\abstract{  We consider locally minimizing measures for the conservative twist maps of the $d$-dimensional annulus or for the Tonelli Hamiltonian flows defined on a cotangent bundle $T^*M$. 
For weakly hyperbolic such measures (i.e. measures with no zero Lyapunov exponents), we prove that the mean distance/angle between the stable and the unstable Oseledet's bundles gives an upper bound of the sum of the positive Lyapunov exponents and a lower bound of the smallest positive Lyapunov exponent.  Some more precise results are proved too.
}
\newpage
\tableofcontents
\newpage
\section{Introduction}
The purpose of this article is to give some relationships between the positive Lyapunov exponents and the angles of the Oseledet's bundles for the invariant minimizing Borel probability measures of the  conservative twisting dynamics.\\

Conservative twisting dynamics are either  what is called a Tonelli Hamiltonian defined on  the cotangent bundle $T^*M$ of a connected and compact manifold $M$ of a twist map of the $d$-dimensional annulus $\A_d=\T^d\times \R^d$. Their two main properties are the following ones:
\begin{enumerate}
\item[--] they twist  the verticals;
\item[--] they are symplectic.
\end{enumerate}
A lot of famous dynamical systems are such conservative twisting dynamics. Let us mention at first all the geodesic flows and mechanical systems (sum of a kinetic energy and a potential one):  they define Tonelli Hamiltonian flows. The twist maps of the two-dimensional annulus were introduced at the end of the nineteenth century by H.~Poincar\'e in the study of the restricted planar circular three body problem (which is a kind of modeling of the system Sun-Earth-Moon). Let us mention too that billiard maps are conservative twist maps and that the Frenkel-Kontorova model can be represented by a conservative twist map (see e.g. \cite{Go1}).\\

To such   dynamics we can associate what is called an {\em action}\footnote{It will be precisely defined later.}, defined via a generating function or a Lagrangian functional.  Two types of action can be defined: either it  is a functional defined along the pieces of orbits or the action of  every  invariant probability measure is defined. The objects of our study are then the locally  minimizing orbits or measures. In the case of Tonelli Hamiltonians, it is well-known that  those orbits (resp. measures) are exactly those that have no conjugate points (see for example \cite{C-I},  \cite{Arna2}). Following \cite{BiaMac}, we will see in subsection \ref{subtwist} that  the locally minimizing orbits of the twist maps of $\A_d$ also have no conjugate points. \\

The following fact is proved in \cite{BiaMac} in the case of symplectic twist maps and in \cite{Arna2}, \cite{C-I} in the case of Tonelli Hamiltonians. \\
If an orbit is locally  minimizing (this means that every piece of this orbit minimizes  locally   the action among the  segments that have  same ends), then there exist  along the orbit two   Lagrangian sub-bundles, invariant under the linearized dynamics and transverse to the vertical bundle, called the Green bundles. These Green bundles enjoy  a lot of nice properties that we will precisely describe  later. Following \cite{Arna2} and \cite{Arna1}, we will denote them by $G_-$ and $G_+$. 

Let us assume now that there exists along such a locally minimizing orbit either an Oseledet's splitting or a partially hyperbolic splitting. We denote the stable, unstable and center bundles corresponding to this splitting by $E^s$, $E^u$ an $E^c$.  It is proved in \cite{Arna4} that: 
$$E^s\subset G_-\subset E^s\oplus E^c\quad{\rm and}\quad E^u\subset G_+\subset E^u\oplus E^c.$$
Hence for a minimizing Borel probability $\mu$, the whole information concerning the positive Lyapunov exponents is contained in the linearized dynamics restricted to the positive Green bundle $G_+$ above the support of $\mu$. Moreover, the angle/distance between the stable and unstable bundles is related to angle/distance between the two Green bundles. \\
Let us recall too that we proved in \cite{Arna4} that for an ergodic locally minimizing measure of a Tonelli Hamiltonian flow, two times the almost everywhere dimension of $G_-\cap G_+$ is equal to the number of zero Lyapunov exponents. The result is  valid for twist maps too.\\

For general dynamical systems, one ``inequality'' between the angles of the Oseledet's splitting and the Lyapunov exponents is well-known; roughly speaking, the smaller the angle/distance between $E^s$ and $E^u$ is, the closer to zero the Lyapunov exponents are. This will be recalled in section \ref{sec3}. In this section too, we will prove  two  exact formulas linking the distance between the two Green bundles and the Lyapunov exponents of the minimizing measures of the conservative dynamics. They are contained in the following theorems.
\begin{thm}\label{thentropy}
Let $\mu$ be a Borel probability measure with no conjugate points that is ergodic for a Tonelli Hamiltonian flow.   If $G_+$ is the graph of $\U$ and $G_-$ the graph of $\stable$, the sum of the positive Lyapunov exponents of $\mu$ is equal to:
$$\Lambda_+(\mu)=\frac{1}{2}\int {\rm tr}(\frac{\partial^2 H}{\partial p^2}(\U-\stable))d\mu.$$
\end{thm}
\noindent Theorem \ref{thentropy} is a slight improvement of a theorem of A.~Freire and R.~Man\'e concerning the geodesic flows that is contained in  \cite{F-M} (see \cite{Fo1} and\cite{C-I} too).  A similar statement was given in the (non published) thesis of G.~Kniepper.\\
Theorem \ref{thentwist} gives a similar statement for the twist maps. In this statement, $G_k(x)=Df^k(f^{-1}(x))V(x)$ is some image of the vertical $V(x)$ that will be precisely defined in section \ref{sec2}.

\begin{thm}\label{thentwist}
Let $f: \A_d\rightarrow \A_d$ be a   twist map and let $\mu$ be a locally minimizing ergodic measure with compact support . Then, if $\Lambda (\mu)$ is the sum of the non-negative exponents of $\mu$,  if $S_-$, $S_+$ designate the symmetric matrices  whose graphs are the  two Green bundles $G_-$ and $G_+$   and $S_k$ designates the symmetric matrix whose graph is $G_k$, we have:

$$\Lambda (\mu)=\frac{1}{2}\int \log\left( \frac{\det\left (S_+(x)-S_{-1}(x)\right)}{\det\left (S_-(x)-S_{-1}(x)\right)}
\right)d\mu (x).$$ 

\end{thm}

For general dynamics, there is no inequality in the other sense. More precisely, the distance between the stable and unstable bundles can be big for measures having Lyapunov exponents that are close to zero. We will see that this phenomenon cannot happen for conservative twisting dynamics. In the following theorems, we denote by $q_+(S)$ the smallest positive eigenvalue of a semi-positive non-zero 
matrix $S$ and we use the same notation as in theorem  \ref{thentropy} for $\U$ and $\esse$.
\begin{thm}\label{th3}
Let $\mu$ be an ergodic     measure with no conjugate points and with at least one non zero Lyapunov exponent for the Tonelli Hamiltonian flow of $H~: T^*M\rightarrow \R$; then its  smallest positive Lyapunov exponent $\lambda(\mu)$ satisfies:
$\lambda(\mu) \geq \frac{1}{2}\int q_+(\frac{\partial^2H}{\partial p^2}).q_+(\U -\stable) d\mu$.
\end{thm}
\noindent Hence, the gap between the two Green bundles gives a lower bound of the smallest positive Lyapunov exponent.\\
For the conservative twist maps, we obtain a similar inequality when all the Lyapunov exponents are non-zero. In this case, the two Green bundles are nothing else but the stable and unstable bundles.

\begin{thm}\label{th4}
Let $f: \A_d\rightarrow \A_d$ be a symplectic twist map and let $\mu$ be a locally minimizing ergodic measure with no zero Lyapunov exponents. We denote the smallest  positive  Lyapunov exponent of $\mu$ by $\lambda (\mu)$ and  an upper bound for $\| s_1-s_{-1}\|$ above ${\rm supp} \mu$ by   $C$. Then we have: 
$$\lambda (\mu)\geq\frac{1}{2}\int \log\left( 1+\frac{1}{C}q_+(\U (x)-\esse (x))\right) d\mu (x).$$
\end{thm}
J.-C.~Yoccoz pointed to me the following illustration of this last result. Let us consider a minimizing fixed point $x_0$ of a two-dimensional twist map $f: \A_1\rightarrow \A_1$. At such a minimizing fixed point, $Df$ has  necessarily two positive eigenvalues denoted by $\lambda$ and $\frac{1}{\lambda}$. Let us denote the matrix of $Df(x_0)$ in the usual coordinates by~: $N=\begin{pmatrix}
a&b\\
c&d\\
\end{pmatrix}$. The twist condition    gives a constant  $\alpha>0$ such that $b\geq \alpha$.  If $Df(x_0)$ is bounded by a constant $C$, this implies that $N$ cannot be to close to the matrix $\begin{pmatrix}
1&0\\
0&1\\
\end{pmatrix}$. Hence $N$ cannot have simultaneously two different eigenvalues close to 1 and a big distance between its eigenspaces. Hence we understand for this example the result contained in the last theorem.\\

Let us comment forward about related results. In \cite{Arna3}, we proved some results concerning the invariant probability measures  of a 2-dimensional twist map whose support is an irrational Aubry-Mather set $\Ac$. We defined at each point $x$ of such an Aubry-Mather set its Bouligand's paratingent cone  $C_x\Ac$, that is a kind of generalized tangent bundle for sets that are not manifolds. We can identify $C_x\Ac$ with the set $\Sc_x\Ac$ of the slopes of its vectors. Then, if $g_-$ and $g_+$ designate the slopes of the two Green bundles, we proved the following inequality~: $g_-(x)\leq \Sc_x\leq g_+(x)$. From that and from theorem \ref{th4}, we deduce~:\\
{\sl The more irregular the Aubry-Mather set is, i.e. the bigger its paratingent cone is, the bigger the Lyapunov exponents are. 

}\medskip

\noindent{\small\sc Acknowledgments}. {\small I am  grateful  to H.~Eliasson and J.-C.~Yoccoz for stimulating discussions and to S.~Crovisier for  pointing to me some improvements of the proofs in section \ref{sub31}.
}

 \section{Some results about the Green bundles}\label{sec2}
\begin{nota}
We assume that $M$ is a compact and connected $d$-dimensional manifold endowed with a fixed Riemannian metric (the associated scalar product is denoted by $(.|.)$). We denote a point of its cotangent bundle $T^*M$ by $(q,p)$ where $p\in T_q^*M$. If $q$ are some (local) coordinates on $M$, then $p$ designate the dual coordinates. This means that if $\eta\in T^*M$ and  $\eta =\sum \eta_idq_i$, then $p_i=\eta_i$.\\
Let us recall that $T^*M$ can be endowed with a 1-form $\lambda$ called the Liouville 1-form, whose expression in all dual coordinates is $\lambda (q,p)=\sum p_idq_i$. Then the canonical symplectic form $\omega$ is defined on $M$ by  $\omega=-d\lambda$. All the dual coordinates are  symplectic for $\omega$.\\
We will denote the usual projection from $T^*M$ to $M$ by $\pi: T^*M\rightarrow M$. For every $x=(q,p)\in T^*M$, we will denote the vertical by $V(x)=\ker (D\pi(x))$. It is a Lagrangian linear subspace of $T_x(T^*M)$.\\
When $M=\T^d$ we will use the global coordinates of $\A_d=\T^d\times \R^d$.
\end{nota}

\subsection{Comparison of two Lagrangian subspaces that are transverse to the vertical}\label{sublag}

Let us recall that a $d$-dimensional subspace $G$ of $T_x(T^*M)$ that is transverse to the vertical $V(x)$ is Lagrangian if and only for every dual linear coordinates $(\delta q, \delta p)$ of $T_x(T^*M)$, then $G$ is the graph of a symmetric matrix in these coordinates.\\

We defined in \cite{Arna2} an order relation for such Lagrangian subspaces of $T_x(T^*M)$ that are transverse to the vertical. The definition is intrinsic and doesn't depend on the chosen dual coordinates,  but let us recall  its interpretation in terms of symmetric matrices: we say that the graph of $\ell_1$ is {\em above} (resp. {\em strictly above}) the graph of $\ell_2$ if the symmetric matrix $\ell_1-\ell_2$ is a positive semi-definite (resp. definite) matrix.\\
We need  to be more precise. Let us recall what we did  in \cite{Arna2}; we associated its {\em height} $Q(S,U)$ to each pair $(S,U)$ of Lagrangian linear subspaces of $T_x(T^*M)$ that are transverse to the vertical. This height is a quadratic form defined on the quotient  linear space $T_x(T^*M)/V(x)$. As this last space is canonically isomorphic to $T_qM$ if $q=\pi (x)$, we modify slightly the set where this quadratic form is defined in comparison with \cite{Arna2}\footnote{We thank F.~Laudenbach for this suggestion.}: \\
if $U$, $S$ are two Lagrangian linear subspaces of $T_x(T^*M)$ that are transverse to the vertical 
$V(x)$, the relative height between $S$ and $U$ is the quadratic form $q(S,U)$ defined on $T_qM$ by the following way: \\

\noindent if $\delta q\in T_qM$, if $\delta x_U\in U$ (resp. $\delta x_S\in S$) is the vector of $U$ (resp. $S$) such that $D\pi(\delta x_U)=\delta q$ (resp. $D\pi (\delta x_S)=\delta q$), then we have: $q(S,U)(\delta q)=\omega (\delta x_S, \delta x_U)$.\\

\noindent Of course, this definition doesn't depend on the dual coordinates that we  choose.  We associate to this bilinear form a unique symmetric operator $s(S, U): T_qM\rightarrow T_qM$ defined by:
$q(S,U)(\delta q_1, \delta q_2)=(s(S,U)\delta q_1| \delta q_2)$.  The operator $s(S,U)$ depends only on the Riemannian product $(.|.)$. Hence, the eigenvalues of $s(S,U)$ are intrinsically defined. We denote them by:
$\lambda_1(S,U)\leq \dots \leq \lambda_d(S,U)$.

\begin{defin}
The quadratic form $q(S,U): T_qM\rightarrow \R$ is called the {\em height } of $U$ above $S$.  The numbers $\lambda_1(S,U)\leq \dots \leq \lambda_d(S,U)$ are the {\em characteristic numbers} of $U$ above $S$.\\
\end{defin}

Let us recall some properties that are proved in \cite{Arna2}. 
\begin{prop} Let $L_1$, $L_2$ and $L_3$ be three Lagrangian subspaces of $T_x(T^*M)$ that are transverse to the vertical. Then:
\begin{enumerate}
\item $\ker q(L_1, L_2)=D\pi (L_1\cap L_2)$;
\item $q(L_1, L_2)=-q(L_2, L_1)$;
\item $q(L_1, L_2)+q(L_2, L_3)=q(L_1, L_3)$.
\end{enumerate}
\end{prop}

\begin{defin}
The {\em distance} between $S$ and $U$ is then \\
 $\displaystyle{\Delta (S, U)= \| q(S,U)\|=\max_{\|\delta q_i\|=1, i=1,2} \omega (\delta x_S^1, \delta x_U^2) }$ where $\delta x_U^i$ (resp. $\delta x_S^i$) designates the element of $U$ (resp. $S$) whose projection on $T_qM$ is $\delta q_i$.

\end{defin}
Let us notice that $\Delta (S,U)$ is not  symplectically invariant.  \medskip

\remk There is a relationship between the distance $\Delta (S, U)$ and the characteristic numbers: $\Delta (S, U)= \max\{|\lambda_1|, |\lambda_d|\}$.
\medskip

\subsection{Tonelli Hamiltonians} We recall some well-known facts concerning Hamiltonian and Lagrangian dynamics (see \cite{Arno}, \cite{Fa2}).\\
\begin{defin} A $C^2$ function $H: T^*M\rightarrow \R$ is called a {\em Tonelli Hamiltonian} if it is:\\
$\bullet$ superlinear in the fiber, i.e. $\forall A\in \R, \exists B\in \R, \forall (q,p)\in T^*M, \| p\|\geq B\Rightarrow  H(q,p) \geq A\| p\|$;\\
$\bullet$ $C^2$-convex in the fiber i.e. for every $(q,p)\in T^*M$, the Hessian $\frac{\partial^2H}{\partial p^2}$ of $H$ in the fiber direction is positive definite as a quadratic form.\\
We denote the Hamiltonian flow of $H$ by $(\varphi_t)$ and the Hamiltonian vector-field by $X_H$. \\
A {\em Lagrangian function} $L: TM\rightarrow \R$ is associated with $H$. It is defined by \\
$\displaystyle{L(q, v)=\max_{p\in T^*_qM} (p.v-H(q,p))}$. 
\end{defin}
Then $L$ is $C^2$-convex and superlinear in the fiber and has the same regularity as $H$. We denote its Euler-Lagrange flow by $(f_t)$. Then $(\varphi_t)$ and $(f_t)$ are conjugated by the Legendre diffeomorphism $\Lc: (q, p)\in T^*M\rightarrow (q, \frac{\partial H}{\partial p}(q, p))\in TM$; more precisely, we have  $\Lc\circ \varphi_t=f_t\circ \Lc$.\\
Let us recall that the orbit of $x\in T^*M$ is $(x_t)=(\varphi_tx)_{t\in\R}$. An {\em  infinitesimal orbit} along the orbit $(x_t)$ is then $(D\varphi_t.\delta x)_{t\in\R}$ where $\delta x\in T_x(T^*M)$. Such an infinitesimal orbit is a solution of the linearized Hamilton equations along the orbit $(x_t)$.

\medskip

The Lagrangian action $A_L(\gamma )$ of a $C^1$ arc $\gamma: [a, b]\rightarrow M$ is defined by:
$$A_L(\gamma )=\int_a^bL(\gamma (s), \dot\gamma (s))ds.$$
A $C^1$ arc  $\gamma_0: [a, b]\rightarrow M$ is {\em minimizing} (resp. {\em locally minimizing})  if for every $C^1$ arc  $\gamma: [a, b]\rightarrow M$ that has the same endpoints as $\gamma_0$, i.e. such that $\gamma_0(a)=\gamma (a)$ and $\gamma_0(b)=\gamma (b)$ (resp. that has the same endpoints as $\gamma_0$, i.e. such that $\gamma_0(a)=\gamma (a)$ and $\gamma_0(b)=\gamma (b)$ and that is sufficiently close to $\gamma_0$ for the $C^1$-topology), we have: $A_L(\gamma_0)\leq A_L(\gamma)$. Such a minimizing (resp.  locally minimizing) arc is   the projection of a unique piece of orbit of the Hamiltonian flow (and then of the Lagrangian flow too). We will say that the corresponding piece of orbit $(\varphi_t(x))_{t\in [a,b]}$ is minimizing (resp. locally minimizing). We say that a complete orbit is minimizing (resp. locally minimizing) if all its restrictions to compact intervals are minimizing (resp. locally minimizing). J.~Mather proved at the end of the 80's (see \cite{mather1}) that there always exist some minimizing orbits. More precisely, he proved the existence of minimizing measures, i.e. Borel probability invariant measures of $T^*M$  whose support is filled with  minimizing orbits.\\

It is well-known that an orbit $(x_t)=(\varphi_t x)$ is locally minimizing if and only if it has no conjugate points. This means that $\forall t\not=u, D\varphi_{t-u}V(x_u)\cap V(x_t)=\{ 0\}$. At every point $y$ of such a minimizing orbit, the family $(G_t(y))=(D\varphi_t.V(\varphi_{-t}y))_{t>0}$ (resp. $(G_{-t}(y))=(D\varphi_{-t}.V(\varphi_{t}y))_{t>0}$ is a decreasing (resp. increasing) family of Lagrangian subspaces that are transverse to the vertical $V(y)$ (see \cite{C-I}, \cite{It1} or \cite{Arna2}) and for every $t>0$, $G_{-t}(y)$ is strictly under $G_t(y)$. Then we define the two Green bundles by
$$G_-(y)=\lim_{t\rightarrow +\infty} G_{-t}(y)\quad{\rm and}\quad G_+(y)=\lim_{t\rightarrow +\infty} G_{t}(y).$$
They are transverse to the vertical, between all the $G_{-t}$ and $G_t$ and $G_+$ is above $G_-$ (see \cite{Arna2} for details).

As at the end of the introduction, let us assume now that there exists along  a locally minimizing orbit either an Oseledet's splitting or a partially hyperbolic splitting. We denote the stable, unstable and center bundles corresponding to this splitting by $E^s$, $E^u$ an $E^c$.  It is proved in \cite{Arna4} and \cite{Arna3} that: 
$$E^s\oplus\R X_H \subset G_-\subset E^s\oplus E^c\quad{\rm and}\quad E^u\oplus \R X_H\subset G_+\subset E^u\oplus E^c.$$

\subsection{Twist maps}\label{subtwist}
The main part of this subsection comes from \cite{BiaMac} (see \cite{Go1} too), even if we changed some proofs. All of what concerns the comparison between the two Green bundles is new. We consider a $C^2$-function $\Phi: \R^d\times \R^d\rightarrow \R$ such that:
\begin{enumerate}
\item $\Phi$ is $\Z^d$-periodic, i.e: $\forall k\in \Z^d, \forall (q, Q)\in\R^d\times\R^d, \Phi(q+k, Q+k)=\Phi(q, Q)$;
\item $\Phi$ satisfies the uniform twist condition, i.e there exists $K>0$ such that:
$$\forall \zeta\in\R^d, \sum_{i, j} \frac{\partial ^2 \Phi(q, Q)}{\partial q_i\partial Q_j}\zeta_i\zeta_j\leq -K\| \zeta\|^2.$$
\end{enumerate}
Then, if we denote the derivative with respect to the $q_i$, $Q_j$ variables by $\Phi_1$ and $\Phi_2$ respectively, the following implicit formula defines a symplectic diffeomorphism $\tilde f$ of $\R^d$:
$$\tilde f(q, p)=(Q, P)\quad{\rm where} \quad P=\Phi_2(q, Q)\quad{\rm and}\quad p=-\Phi_1(q,Q).$$
We say then that $\Phi$ is a {\em generating function} for $\tilde f$. We associate a formal function defined on $(\R^d)^\Z$ to $\Phi$:
$$\Ac((q_n)_{n\in\Z})=\sum_{n=-\infty}^{+\infty} \Phi(q_n, q_{n+1}).$$
Even if this function is not well-defined, its critical points  are well-defined, they satisfy the equations:
$$\forall n\in\Z, \Phi_2(q_{n-1}, q_n)+\Phi_1(q_n, q_{n+1})=0.$$
We can denote the partial actions $\Ac_{M, N}$ for $M\leq N$ by:
$$\Ac_{M, N}((q_n)_{M\leq n\leq N})=\sum_{n=M}^{N-1} \Phi(q_n, q_{n+1}).$$
Then $(q_n)_{M\leq n\leq N}$ is a critical point of $\Ac_{M, N}$ restricted to the set of the finite sequences that have the same endpoints as $(q_n)_{M\leq n\leq N}$ if and only if it is a the projection of a finite piece of orbit $(q_n, p_n)_{M\leq n\leq N}$ for $\tilde f$.  In this case, we have:
\begin{enumerate}
\item[$\bullet$] $p_M=-\Phi_1(q_M, q_{M+1})$; $p_N=\Phi_2(q_{N-1}, q_N)$;
\item[$\bullet$] $\forall n\in [M+1, N-1], p_n=\Phi_2(q_{n-1}, q_n)=-\Phi_1(q_n, q_{n+1})$.
\end{enumerate}
We say that $(q_n)_{M\leq n\leq N}$ is {\em minimizing} (resp. {\em locally minimizing}) if it is minimizing (resp. locally minimizing) among all the segments that have the same endpoints. Then the corresponding piece of orbit $(q_n, p_n)_{M\leq n\leq N}$ is said to be minimizing (resp. locally minimizing) too. We say that $(q_n)_{n\in\Z}$ or $(q_n, p_n)_{n\in\Z}$ is minimizing (resp.  locally minimizing) if all its restrictions to segments are minimizing (resp. locally minimizing).

If now $(x_n)=(q_n, p_n)\in (\A_d)^\Z$ is an orbit for $f$, we say that it is {\em minimizing} (resp. {\em locally minimizing}) if its lifted orbit $(\tilde q_n, p_n)$ for $\tilde f$ is minimizing. Moreover, we will denote the partial action of the lift by: 
$$\Phi_{N, M}((q_n))=\Ac_{N, M}((\tilde q_n)).$$\\

Let us now fix an orbit $(x_n)=(q_n, p_n)$ for $f$. We call an {\em infinitesimal orbit} along $(x_n)$ a sequence $(Df^n(x_0)\delta x)_{n\in\Z}$, i.e. an infinitesimal orbit is an orbit for the derivative of $f$. The projection of an infinitesimal orbit is called a {\em Jacobi field}. Then $(\zeta_n)$ is a Jacobi field if and only if we have:
$$\forall n\in\Z, {}^tb_{n-1}\zeta_{n-1}+a_n\zeta_n+b_n\zeta_{n+1}=0;$$
where $b_n=\Phi_{12}(q_n, q_{n+1})$ and $a_n=\Phi_{11}(q_n, q_{n+1})+\Phi_{22}(q_{n-1}, q_n).$\\
The Hessian of $\Phi_{M, N}$ is:
$$D^2\Phi_{M, N}((x_n))=\begin{pmatrix}
a_M& b_M&0&\dots &\dots&\dots&0\\
{}^tb_M&a_{M+1}&b_{M+1}&0&\dots&\dots&0\\
\dots&\dots&\dots&\dots&\dots&\dots&\dots\\
\dots&\dots&\dots&\dots&\dots&a_{N-1}&b_{N-1}\\
0&\dots&\dots&\dots&0&{}^tb_{N-1}&a_N\\
\end{pmatrix}.
$$
The kernel of this Hessian is made with the Jacobi fields $(\zeta_n)_{M\leq n\leq N}$ such that $\zeta_{M-1}=\zeta_{N+1}=0$.\\
If we assume that $(x_n)$ is locally minimizing, then all the Hessians $D^2\Phi_{M, N}((x_n))$ are a priori positive semi-definite.  Following \cite{BiaMac}, let us prove that these Hessians are in fact positive definite.
\begin{prop} {\bf (Bialy-MacKay, \cite{BiaMac})}
If the orbit $(x_n)$ of $f$ is locally minimizing, then all the Hessians $D^2\Phi_{M, N} ((x_n))$ are positive definite and then the orbit has no conjugate vectors.
\end{prop}
\demo
If not, there exist  $M\leq N$   and a Jacobi field $(\zeta_n)_{n\in\Z}$ that is different from $(0)$ but such that $\zeta_{M-1}=\zeta_{N+1}=0$.  In other words, this Jacobi field has what is usually called {\em conjugate vectors}. In this case, $(0, 0, \zeta_M, \zeta_{M+1}, \dots , \zeta_{N-1}, \zeta_N, 0, 0)$ is in the isotropic cone of $D^2\Phi_{M-2, N+2}((x_n))$ but not in its kernel (because it is not a Jacobi field); this contradicts the fact that the kernel is equal to the isotropic cone (because this Hessian is positive semi-definite).\enddemo
\noindent Hence the Jacobi fields along any locally minimizing orbit have no conjugate vectors. This implies  that for any $k\in\Z^*$ and any $n\in\Z$, $G_k(x_{n+k})=Df^k(x_n).V(x_n)$ is transverse to $V(x_{n+k})=V(f^kx_n)$. \begin{prop}{\bf (Bialy-MacKay, \cite{BiaMac})}
Let $(x_k)$ be a locally minimizing orbit. Then, for all $k\geq 1$, we have along this orbit:
\begin{enumerate}
\item[$\bullet$] $G_{-1}$ is strictly under $G_{k}$ and $G_{-k}$ is strictly under $G_1$;
\item[$\bullet$] $G_{k+1}$ is strictly under $G_k$ and $G_{-k}$ is strictly under $G_{-(k+1)}$.
\end{enumerate}
\end{prop}
This result is proved in  \cite{BiaMac}, but we give a slightly different proof.\\

We deduce that  $(G_k)_{k\geq 1}$ is a decreasing sequence of Lagrangian subspaces that are all above $G_{-1}$, hence we can define $\displaystyle{G_-=\lim_{k\rightarrow +\infty}G_k}$. Similarly, $(G_{-k})_{k\geq 1}$ is an increasing sequence of Lagrangian subspaces that are all under $G_1$, hence we can define $G_+$ by taking the limit.

\begin{defin}
If the orbit of $x$ is locally minimizing, the two {\em Green bundles} at $x$ are the two Lagrangian subspaces of $T_x(T^*M)$ that are transverse to the vertical and defined by:
$$G_-(x)=\lim_{k\rightarrow +\infty} G_{-k}(x)\quad{\rm and}\quad G_+(x)=\lim_{k\rightarrow +\infty} G_k(x).$$
\end{defin}
\demo
We denote the symmetric matrix whose graph is $G_k(x_{n+k})$  by  $S_k(x_{n+k})$.
\\
Let us notice that: $$ Df(x_n)=\begin{pmatrix}
-b_n^{-1} \Phi_{11}(q_n, q_{n+1})&-b_n^{-1}\\
{}^tb_n-\Phi_{22}(q_n, q_{n+1})b_n^{-1}\Phi_{1 1}(q_n, q_{n+1}) & -\Phi_{2 2} (q_n, q_{n+1})b_n^{-1}\\
\end{pmatrix}.$$
We deduce that: $G_1(x_{n+1})={\rm graph} (\Phi_{22}(q_n, q_{n+1}))$,  $G_{-1}(x_n)={\rm graph} (-\Phi_{1 1}(q_n, q_{n+1}))$ and then   $S_1(x_n)=\Phi_{2 2}(q_{n-1}, q_n)$, $S_{-1} (x_n)=-\Phi_{1 1}(q_n, q_{n+1})$. Hence:  $a_n=S_{1}(x_n)-S_{-1}(x_n)$ is the matrix of the relative height between   $G_{-1}(x_n)$ and $G_1(x_n)$   (see subsection \ref{sublag} for definition). Hence $G_1$ is strictly above $G_{-1}$. 

Let us prove that: $\forall k\geq 1, S_k(x_n)-S_{-1}(x_n)>0$.  If not, there exists $k\geq 2$ and $\eta\not=0$ such that: ${}^t\eta (S_k(x_n)-S_{-1}(x_n))\eta\leq0$. Then we consider the piece of infinitesimal orbit $\left(Df^{-j}\begin{pmatrix}\eta\\ S_k(x_n)\eta\
\end{pmatrix}\right)_{0\leq j\leq  k-1}$ and the Jacobi field that is the projection of this infinitesimal orbit: $\zeta_i=D\pi\circ Df^{i-n} \begin{pmatrix}\eta\\ S_k(x_n)\eta\
\end{pmatrix}$ for $n-k+1\leq i\leq  n$. Let us compute $D^2\Phi_{n-k+1, n}(x)\zeta=\Delta$. \begin{enumerate}
\item  as $Df^{-k}G_k(x_n)=V(x_{n-k})$, we have: $\Delta_{n-k+1}=a_{n-k+1}\zeta_{n-k+1}+b_{n-k+1}\zeta_{n-k+2}=-{}^tb_{n-k}\zeta_{n-k}=0$;
\item as we have  a Jacobi field, for $n-k+1\leq i\leq n-2$, we have:    $\Delta_{i+1}= {}^tb_i\zeta_i+a_{i+1}\zeta_{i+1}+b_{i+1}\zeta_{i+2}=0$;
\item $\Delta_n={}^tb_{n-1}\zeta_{n-1}+a_n\zeta_n=-b_nD\pi\circ Df \begin{pmatrix}\eta\\ S_k(x_n)\eta\
\end{pmatrix}$\\
$=-b_n\left( -b_n^{-1}(\Phi_{11}(q_n, q_{n+1})-S_k(x_n))
\right)\eta=-(S_{-1}(x_n)-S_k(x_n))\eta$.
\end{enumerate}
We deduce that $D^2\Phi_{n-k+1, n}(x)(\zeta, \zeta)={}^t\Delta.\zeta={}^t\eta (S_k(x_n)-S_{-1}(x_n))\eta\leq 0$. This contradicts the fact that the Hessian is positive definite. Hence we have proved that for all positive $k$, $G_k$ is strictly above $G_{-1}$.

Moreover,    $G_{k+1}(x_{n+1})$ is represented by: \\
$Df(x_n)\begin{pmatrix} {\bf 1}\\ S_k(x_n)\\ \end{pmatrix}=\begin{pmatrix}
-b_n^{-1}(\Phi_{11}(q_n, q_{n+1})+S_k(x_n))\\
{}^tb_n-\Phi_{22}(q_n, q_{n+1})b_n^{-1}(\Phi_{11}(q_n, q_{n+1})+S_k(x_n))\\
\end{pmatrix}.$\\
This means: $S_{k+1}(x_{n+1})=-{}^tb_n(\Phi_{11}(q_n , q_{n+1})+S_k(x_n))^{-1}b_n+\Phi_{22}(q_n , q_{n+1})$ and then:\\
$(S_{k+1}-S_{-1})(x_{n+1})= a_{n+1}-{}^tb_n((S_k-S_{-1})(x_n))^{-1}b_n$ i.e.: \\
$(S_{k+1}-S_{-1})(x_{n+1})=(S_1-S_{-1})(x_{n+1})-{}^tb_n((S_k-S_{-1})(x_n))^{-1}b_n$.\\
In particular, we have: $(S_{2}-S_{-1})(x_{n+1})=(S_1-S_{-1})(x_{n+1})-{}^tb_na_n^{-1}b_n$ then $S_2<S_1$.\\
We can subtract for any $k\geq 2$: 
$$(S_{k+1}-S_k)(x_{n+1})={}^tb_n\left((S_{k-1}-S_{-1})(x_n))^{-1}-(S_k-S_{-1})(x_n)^{-1}\right) b_n.$$
 We have proved that for all positive $k$, $G_k$ is strictly above $G_{-1}$. We deduce that $(G_k(x_n))_{k\geq 1}$ is a strictly decreasing sequence of Lagrangians subspaces. Because all these subspaces are above $G_{-1}(x_n)$, they converge to a Lagrangian subspace $G_+$ that is transverse to the vertical. In the same way, we obtain that $(G_{-k}(x_n)_{k\geq 0}$ is an increasing sequence of Lagrangian subspaces that are bounded from above by $G_1$, hence they converge to a Lagrangian subspace $G_-$ that is transverse to the vertical.\enddemo
 \begin{prop}
 Let $x\in T^*M$ whose orbit is locally minimizing. Then for all $n, k\geq 1$, $G_{-k}(x)$ is strictly under $G_n(x)$.\\
 Hence $G_-$ is under $G_+$. \end{prop}
 \demo We denote $f^m(x)$ by $x_m$.  Let us   prove that for all $n, k\geq 1$, and all $m\in\Z$, then $G_n(x_m)$ is above $G_{-k}(x_m)$. We have proved this result for $k=1$ or $n=1$, then we assume that $n, k\geq 2$.\\
 We recall that if $F_1$, $F_2$ are two transverse Lagrangian subspaces of a symplectic space whose dimension is denoted by $2d$, then the set $\Tc (F_1, F_2)$ of the Lagrangian subspaces that are transverse to both $L_1$ and $L_2$ has exactly $d+1$ connected components: it depends on the signature of a certain quadratic form. Let us consider the connected component $\Cc$ of $\Tc(G_{k-1}(x_{k-1+m}), G_{k-1+n}(x_{k-1+m}))$ that contains  $G_{k+n}(x_{k-1+m})$;  we have proved that $G_{k+n}(x_{k-1+m})$ and $G_{-1}(x_{k-1+m})$ are under $G_{k-1}(x_{k-1+m})$ and $ G_{k-1+n}(x_{k-1+m})$, hence they are in the same connected component $\Cc$ of  $\Tc(G_{k-1}(x_{k-1+m}), G_{k-1+n}(x_{k-1+m}))$ and their images by $\left(Df^{k-1}(x_m)\right)^{-1}$, that are $G_{n+1}(x_m)$ and $G_{-k}(x_m)$, are in the same connected component of $$\left(Df^{k-1}(x_m)\right)^{-1}\left( \Tc(G_{k-1}(x_{k-1+m}), G_{k-1+n}(x_{k-1+m}))\right).$$This last set is equal to: \\
 $\Tc(  \left(Df^{k-1}(x_m)\right)^{-1}(G_{k-1}(x_{k-1+m}),  \left(Df^{k-1}(x_m)\right)^{-1}(G_{k-1+n}(x_{k-1+m})))
 $\\
 $=\Tc(V(x_m), G_{n}(x_m)).$\\
 We have proved that $G_{n+1}(x_m)$ is under  $G_{n} (x_m)$. As $G_{n+1}(x_m)$ and $G_{-k}(x_m)$ are in the same connected component of $\Tc(V(x_m), G_{n}(x_m))$, this implies that $G_{-k}(x_m)$ is under $G_{n}(x_m)$. \\
We deduce that $G_-$ is under $G_+$.\enddemo

\section{Sum of the positive Lyapunov exponents and upper bounds}\label{sec3}
Before explaining which results we obtain for the twisting dynamics, we have to explain that some results are true for general dynamics (not necessarily twisting) and explain the difference with our results.
\subsection{ Some general results}\label{sub31}
In this section, we review some more or less well-known results concerning the link between the Lyapunov exponents and the   distance between the Oseledet's bundles. Because we didn't find any precise reference and because the proofs are rather short, we give here a proof a these results. A good reference for Lyapunov exponents is \cite{Led}.

We work on a  manifold $N$ (not necessarily compact) and we ask ourselves the following question.

\noindent{\bf Question.}  If the Oseledet's splitting of an invariant measure of a $C^1$-diffeomorphism is such that $E^s$ and $E^u$ are close to each other (in a sense we have to specify), are the Lyapunov exponents all close to $0$?

Let us explain that the answer is yes if $E^s$ and $E^u$ are $1$-dimensional.

\begin{nota}
If $E$, $F$ are two linear subspaces of $T_xN$ that are $d$-dimensional with $d\geq 1$, the distance between $E$ and $F$ is:  
$${\rm dist} (E,F)=  \inf_{(e_i), (f_i)}\max\{ \| e_1-f_1\|, \dots , \| e_d-f_d\|\}  $$
where the infimum is taken over all the orthonormal basis $(e_i)$ of $E$, $(f_i)$ of $F$.
\end{nota}
\begin{prop}
Let $K$ be a compact subset of $N$ and let $C>0$   be a real number. Then, for any $f\in {\rm Diff}^1(M)$ so that $\max \{ \| Df_{|K}\|, \| Df^{-1}_{|K}\|\}\leq C$, if $f$ has an invariant ergodic measure $\mu$ with support in $K$ such that the Oseledet's stable  and unstable bundles $E^s$ and $E^u$ of $\mu$ are one dimensional, if we denote by $\lambda_u$   the positive  Lyapunov exponent  and by $\lambda_s$ the   negative one, then: 
$$0< \lambda_u-\lambda_s\leq \log\left( 1+C^2\int {\rm dist}(E^u, E^s)d\mu\right).$$
\end{prop}
\demo
We denote ${\rm dist}(E^u( x), E^s( x))$ by $\alpha (x)$.We choose $x\in {\rm supp}\mu$ where $E^s$ and $E^u$ are defined and we choose $v\in E^u(x)\backslash \{ 0\} $. Then there exists $p_x(v)\in E^s(x)$ such that $\| p_x(v)\|=\| v\|$ and $\| p_x(v)-v\|\leq \alpha (x)\| v\|$. Then,  we have:
$$\begin{matrix}  \| Df_{|E^u(x)}\|.\| v\|= \| Df(x)v\|& \leq& \| Df(x)p_x(v)\| +\| Df(x)\| .\| p_x(v)-v\|\\
&\leq&{ \| Df_{|E^s}(x)\| \left( 1+\alpha (x)\frac{\| Df(x)\|}{\|Df_{|E^s}(x)\|}\right)\| v\|.\hfill}\\
\end{matrix}$$
We deduce:
$$\begin{matrix}
\lambda_u-\lambda_s&=&\int\log\| Df_{|E^u}\| d\mu-\int\log\| Df_{|E^s}\| d\mu\\
&\leq& \int\log(1+C^2\alpha(x))d\mu(x)\leq \log\left( 1+C^2\int\alpha (x)d\mu (x)\right) \\
\end{matrix}$$
by Jensen inequality.\enddemo
In the higher dimension cases, we obtain a slightly less good estimation.
\begin{prop}
Let $K$ be a compact subset of $N$,  let $C>0$   be a real number. Then, for any $f\in {\rm Diff}^1(M)$ so that $\max \{ \| Df_{|K}\|, \| Df^{-1}_{|K}\|\}\leq C$, if $f$ has an invariant ergodic measure $\mu$ with support in $K$ such that the Oseledet's stable  and unstable bundles $E^s$ and $E^u$ of $\mu$ have the same dimension $d$  , if we denote by $\Lambda_u$   the sum of the positive   Lyapunov exponents and by $\Lambda_s$ the sum of the negative Lyapunov exponents, then: 
$$0< \Lambda_u-\Lambda_s\leq  d\log\left(1+(C^2+1)\int {\rm dist}(E^u, E^s)d\mu \right).$$
\end{prop}
\demo
We denote ${\rm dist}(E^u( x), E^s( x))$ by $\alpha (x)$. At all the points where $E^s$ and $E^u$ are defined, we choose an orthonormal basis $(e_1, \dots , e_d)$ of $E^s$ in a measurable way, and  an   orthogonal basis $(f_1, \dots , f^d)$ of $E^u$ that depends measurably on the considered point and is such that:
$${\rm dist} (E^s,E^u)=\max \{ \| e_1-f_1\|, \dots , \| e_d-f_d\|\}. $$
Then we denote by $P_x: E^u(x)\rightarrow E^s(x)$ the linear map such that $P_x(f_i(x))=e_i(x)$; then, each $P_x$ is an isometry. Moreover, because $\| P_x(f_i(x))-f_i(x)\|=\| f_i(x)-e_i(x)\|\leq \alpha (x)$, we deduce: $\forall v\in E^u(x), \|  P_x(v)-v\|\leq d\alpha (x)\| v\|$ and $\forall v\in E^s(x), \|  (P_x)^{-1}(v)-v\|\leq d\alpha (x)\| v\|$.
We have (we compute the determinant in the previous basis):
$$\begin{matrix}
\Lambda_u-\Lambda_s&=&\int\left(\log(|\det Df_{|E^u(x)}|)-\log(|\det Df_{|E^s(x)}|)\right) d\mu (x)\\
&=& \int\log\left | \det(Df_{|E^u(x)}(P_{x})^{-1}(Df_{|E^s(x)})^{-1}P_{f(x)})\right |Êd\mu (x).
\end{matrix}$$
Let us consider $v\in E^u(x)$. Then: $$
  \begin{matrix}\|  Df (x) (P_{x})^{-1}&(Df (x))^{-1}P_{f(x)}v-v\| \hfill\\
 \hfill  \leq& \left\| Df(x)\left((P_{x})^{-1}(Df (x))^{-1}P_{f(x)}v-(Df (x))^{-1}P_{f(x)}v)
\right) \right\|+\left\| P_{f(x)}v-v\right\|\\
 \hfill\leq&Cd\alpha (x)\| (Df (x))^{-1}P_{f(x)}v\| +d\alpha(f(x))\| v\|\leq d(C^2\alpha (x)+\alpha (f(x)))\| v\|.\\
\end{matrix}
$$
Hence we have: $\| Df (x) (P_{x})^{-1}(Df (x))^{-1}P_{f(x)}-{\rm Id}_{E^u(x)}\|\leq d(C^2\alpha (x)+\alpha (f(x)))$. We deduce: $$\left|\det\left( Df_{|E^u(x)}(P_{x})^{-1}(Df_{|E^s(x)})^{-1}P_{f(x)}
\right)\right|\leq \left( 1+d(C^2\alpha (x)+\alpha(f(x)))\right) ^d $$
and then:
$$\begin{matrix}
\Lambda_u-\Lambda_s&\leq& d\int\log\left( 1+d(C^2\alpha (x)+\alpha(f(x)))\right) d\mu(x)\\
&\leq& d\log\left(1+(C^2+1)\int\alpha(x)d\mu(x)\right).
\end{matrix}
$$\enddemo

Hence the fact that we will obtain a upper bound of the sum of the positive Lyapunov exponents that depends on the distance between the two Green bundles in the case of the twisting dynamics is not surprising.  The results contained in section \ref{secmin}, that give lower bounds  that are specific to the twisting dynamics, are more surprising. \\
What is more interesting in this section is that we obtain some exact formula for the sum of the positive Lyapunov exponents.

\subsection{Tonelli Hamiltonians}\label{sec32}
Using a Riemannian metric on $M$, we define the horizontal subspace $\Hc$ as the kernel of the connection map. Then, for every Lagrangian subspace $\Gc$ of $T_x(T^*M)$, there exists a linear map $G: \Hc(x)\rightarrow V(x)$ whose graph is $\Gc$. That is the meaning of {\em graph} in the following theorem.\\
\remk If $K$ is an invariant compact and locally minimizing subset of $T^*M$ (for example the support of a  locally minimizing ergodic measure), we have:
$$\forall x\in K, G_{-1}(x)\leq G_-(x)\leq G_+(x)\leq G_1(x)$$
where ``$\leq$'' designates the relation ''to be below'' for the Lagrangian subspaces that are transverse to the vertical.\\
Hence $G_-$ and $G_+$ are uniformly bounded on $K$ and $D\pi: G_\pm(x)\rightarrow T_xM$ is uniformly bilipschitz.

In the case of ergodic measures of a geodesic flow with support filled by locally minimizing orbits, i.e. in the case of measures with no conjugate points, A.~Freire and R.~Man\'e proved in \cite{F-M} a nice formula for the sum of the positive exponents (see \cite{Fo1} and \cite{C-I} too). A slight improvement of this formula gives:

\begin{theo}{\bf \ref{thentropy}}
Let $\mu$ be a Borel probability measure with no conjugate points that is ergodic for a Tonelli Hamiltonian flow.   If $G_+$ is the graph of $\U$ and $G_-$ the graph of $\stable$, the sum of the positive Lyapunov exponents of $\mu$ is equal to:
$$\Lambda_+(\mu)=\frac{1}{2}\int {\rm tr}(\frac{\partial^2 H}{\partial p^2}(\U-\stable))d\mu.$$
\end{theo}
Hence, we see that the more distant the Green bundles are, the greater  the sum of the positive Lyapunov exponents is. This gives an upper bound to the positive Lyapunov exponents.  \\
\demo
A consequence of the linearized Hamilton equations  is that if the graph ${\cal G}$ of a symmetric matrix $G$ is invariant by the linearized flow, then any infinitesimal orbit $(\delta q, G\delta q)$ satisfies the following equation: $\delta\dot q=(\frac{\partial^2H}{\partial p^2}G+\frac{\partial^2 H}{\partial q\partial p})\delta q$.\\
Hence, we have: $\frac{d{}}{dt} \det(D\pi\circ D\varphi_{t|{\cal G}})={\rm tr}(\frac{\partial^2H}{\partial p^2}G+\frac{\partial^2 H}{\partial q\partial p})\det(D\pi\circ D\varphi_{t|{\cal G}})$; we deduce: \\
$$\begin{matrix} \frac{1}{T}\log \det  (D\pi\circ   D\varphi_{T|{\cal G}}&(q,p)) \hfill \\
&\hfill =\frac{1}{T}\log\det (D\pi(q,p)_{|{\cal G}})+\frac{1}{T}\int_0^T {\rm tr}(\frac{\partial^2H}{\partial p^2}G+\frac{\partial^2 H}{\partial q\partial p})(\varphi_t(q,p))dt.\end{matrix}$$
Via ergodic Birkhoff's theorem, we deduce for $(q,p)$ generic that: 
$$\liminf_{T\rightarrow +\infty}\frac{1}{T}\log\det (D\pi\circ D\varphi_{T|{\cal G}}(q,p))=\int {\rm tr}(\frac{\partial^2H}{\partial p^2}G+\frac{\partial^2 H}{\partial q\partial p})d\mu.$$
We have noticed that $D\pi_{G_\pm}$ is uniformly bilipschitz above ${\rm supp}\mu$, hence we can remove $D\pi$ in the previous formula when ${\cal G}$ is one of the two Green bundles.\\
Moreover, we know that $E^s\subset G_-\subset E^{s\bot}=E^c\oplus E^s$ and that $E^u\subset G_+\subset E^{u\bot}=E^c\oplus E^u$. Hence, the sum of the Lyapunov exponents of the restricted cocycle $(D\varphi_{t|G_+})$ is exactly $\Lambda_+(\mu)$  and the sum of the Lyapunov exponents of the restricted cocycle $(D\varphi_{t|G_-})$ is $\Lambda_-(\mu)=-\Lambda_+(\mu)$. Then  we have:
$$\Lambda_+(\mu)=\int{\rm tr}(\frac{\partial^2H}{\partial p^2}\U+\frac{\partial^2 H}{\partial q\partial p})d\mu\quad {\rm and}\quad -\Lambda_+(\mu)=\int{\rm tr}(\frac{\partial^2H}{\partial p^2}\stable+\frac{\partial^2 H}{\partial q\partial p})d\mu .$$
We obtain the conclusion by subtracting   the two equalities.
\enddemo

\subsection {Twist maps}

\begin{theo}{\bf\ref{thentwist}}
Let $f: \A_d\rightarrow \A_d$ be a   twist map and let $\mu$ be a locally minimizing ergodic measure with compact support. Then, if $\Lambda (\mu)$ is the sum of the non-negative exponents of $\mu$,  if $S_-$, $S_+$ designate the symmetric matrices  whose graphs are the  two Green bundles $G_-$ and $G_+$ and $S_k$ designates the symmetric matrix whose graph is $G_k$, we have:

$$\Lambda (\mu)=\frac{1}{2}\int \log\left( \frac{\det\left (S_+(x)-S_{-1}(x)\right)}{\det\left (S_-(x)-S_{-1}(x)\right)}
\right)d\mu (x).$$ 

\end{theo}

In this case again, we see that the closer the two Green bundles are to each other, the closer to $0$  the Lyapunov exponents are.

\demo We use coordinates such that $G_+$ becomes the horizontal bundle, i.e. we use the change of symplectic coordinates whose matrix is $\begin{pmatrix} 1&S_+\\ 0&1\\ \end{pmatrix}$. This change of coordinates is not continuous,  but it is uniformly bounded because $S_{1}\leq S_+\leq S_1$, and $S_{-1}$ and $S_1$ vary continuously in the compact set ${\rm supp}\mu$.
The matrix of $Df$ at $x$  is then: $$M=\begin{pmatrix}
B_1(x)(S_+(x)-S_{-1}(x))&B_1(x)\\
0&B_1(x)(S_1(fx)-S_+(fx))\\
\end{pmatrix}.$$
We know that $E^u\subset G_+\subset E^u\oplus E^c$, hence along $G_+$ we see   the non-negative Lyapunov exponents. Then we have:
$$\Lambda(\mu)=\int\log\left| \det Df_{|G_+}\right| d\mu=\int\log\left| \det B_1(x)(S_+(x)-S_{-1}(x))\right| d\mu .$$
In the same way, we have: 
$$-\Lambda(\mu)=\int\log\left| \det Df_{|G_-}\right| d\mu=\int\log\left| \det B_1(x)(S_-(x)-S_{-1}(x))\right| d\mu .$$

By subtracting these two equalities, we obtain the equality of the theorem.

\enddemo

\section{Lower bounds for the positive Lyapunov exponents}\label{secmin}

Here we prove results that are specific to the twisting dynamics.

\subsection{Tonelli Hamiltonians}
 \begin{lemma}\label{lem9} Let $H: T^*M\rightarrow \R$ be a Tonelli Hamiltonian. Let $(x_t)$ be a locally minimizing orbit and let $U$ and $S$ be two Lagrangian bundles along this orbit that are invariant by the linearized Hamilton flow and transverse to the vertical. Let $\delta x_U\in U$ be an infinitesimal orbit contained in the bundle $U$ and let us denote by $\delta x_S$ the unique vector of $S$ such that $\delta x_U-\delta x_S\in V$ (hence $\delta x_S$ is {\em not} an infinitesimal orbit). Then:
$$\frac{d{}}{dt}(\omega(x_t)(\delta x_S(t), \delta x_U(t))=(\delta x_U(t)-\delta x_S(t))  H_{pp} (x_t)(  \delta x_U(t)-\delta x_S(t))\geq 0.$$

\end{lemma}
\remk \quad  Let us notice that $\omega (x_t)(\delta x_S(t), \delta x_U(t))$ is nothing else  but the relative height $Q(S,U)$ of $U$ above $S$ at the vector $\delta q_U=D\pi.\delta x_U=D\pi. \delta x_S$.

\demo As the result that we want to prove is local, we can assume that we are in the domain of a dual chart and express all the things in the corresponding dual linearized coordinates.\\
We consider an invariant Lagrangian  linear bundle $G$ that is transverse to  the vertical along the orbit of $x=(q,p)$. We denote the symmetric matrix whose graph is $G$  by $G$ again. An infinitesimal orbit contained in this bundle satisfies: $\delta p = G \delta q$. We deduce from the linearized Hamilton equations (if we are along the orbit $(q(t), p(t))=x(t)$, $\dot G$ designates $\frac{d{ }}{dt} (G(x(t)))$) that:
$$\delta \dot q = (H_{qp} +H_{pp} G)\delta q; \quad \delta \dot p = (\dot G +GH_{qp}+GH_{pp}G)\delta q=-(H_{qq}  +H_{pq}G)\delta q.$$
We deduce from these equations the classical Ricatti equation (it is  given for example in \cite{C-I} for Tonelli Hamiltonians, but the reader can find the initial and simpler Ricatti equation given by Green in the case of geodesic flows in \cite{Gr}):
$$\dot G +GH_{pp}G+GH_{qp}+H_{pq}G+H_{pp}=0.$$
Let us assume now that the graphs of the symmetric matrices $\U$ and $\esse$ are invariant by the linearized flow  along the same orbit. We denote by $(\delta q_U, \U \delta q_U)$ an infinitesimal orbit that is contained in the graph of $\U$. Then we have: 
$$\frac{d{}}{dt}(\delta q_U(\U-\esse) \delta q_U)=2\delta q_U(\U-\esse ) \delta\dot q_U +\delta q_U(\dot\U -\dot \esse )\delta q_U
$$
$$
=2\delta q_U(\U-\esse )(H_{qp} +H_{pp} \U )\delta  q_U +\delta q_U( \esse H_{pp}\esse - \U H_{pp} \U +\esse H_{qp}+H_{pq}\esse -\U H_{qp}-H_{pq}\U)\delta q_U
$$
$$
=\delta q_U(\U H_{qp} -\esse H_{qp}+\U H_{pp}\U -2\esse H_{pp} \U +\esse H_{pp}\esse  +H_{pq}\esse -H_{pq}\U )\delta q_U
$$
$$
=\delta q_U (\U -\esse)H_{pp} (\U-\esse) \delta q_U \geq 0.
$$
To finish the proof, we just need to notice that in coordinates:
$$\omega(\delta x_S, \delta x_U)=\omega (\delta x_U, \delta x_U-\delta x_S)=(\delta q_U, \U \delta q_U)\begin{pmatrix} 0& {\bf 1}\\
-{\bf 1} & 0\\
\end{pmatrix} \begin{pmatrix} 0\\ (\U-\esse)\delta q_U\\
\end{pmatrix}=\delta q_U (\U-\esse)\delta q_U
$$
\enddemo

\begin{nota}
If $S$ is a positive semi-definite matrix that is not the null matrix, then $q_+(S)$ is its smallest positive eigenvalue.
\end{nota}

\begin{theo}{\bf \ref{th3}}
Let $\mu$ be an ergodic     measure with no conjugate points and with at least one non zero Lyapunov exponent; then its  smallest positive Lyapunov exponent $\lambda(\mu)$ satisfies:
$\lambda(\mu) \geq \frac{1}{2}\int m(\frac{\partial^2H}{\partial p^2}).q_+(\U -\stable) d\mu$.
\end{theo}

Hence, the gap between the two Green bundles gives a lower bound of the smallest positive Lyapunov exponent. It is not surprising that when $E^s$ and $E^u$ collapse, the Lyapunov exponents are $0$. What is more surprising and specific to the case of Tonelli Hamiltonians is the fact that the bigger the gap between $E^s$ and $E^u$ is,  the greater the Lyapunov exponents are: in general, along a hyperbolic orbit, you may have a big angle between the Oseledet's bundles and some very small Lyapunov exponents.

\noindent\begin{proof} Let   $\mu$ be an ergodic Borel probability measure with no conjugate points; its support $K$ is  compact and then, by the first remark of section \ref{sec32}, there exists a constant $C>0$ such that $\U$ and $\esse$ are bounded by $C$ above $K$. We choose a point $(q, p)$ that is generic for $\mu$ and $(\delta q, \U \delta q)$ in the   Oseledet's bundle corresponding to the smallest positive Lyapunov exponent $\lambda(\mu)$ of $\mu$. 
Using the linearized Hamilton equations (see lemma \ref{lem9}), we obtain:
 $$\frac{d{}}{dt}((\delta q(\U -\stable)\delta q)=
\delta q(\U -\stable)\frac{\partial^2H}{\partial p^2}(q_t, p_t)(\U-\stable)\delta q.$$
 Let us notice that $(\U-\stable )^\frac{1}{2}\delta q$ is contained in  the orthogonal space to the kernel of $\U-\stable$. Hence:
  $$\frac{d{}}{dt}((\delta q(\U -\stable)\delta q)\geq m(\frac{\partial^2H}{\partial p^2})q_+(\U-\stable ) \delta q(\U -\stable)\delta q.$$
Moreover $\delta q\notin \ker (\U-\esse)$ because $(\delta q,\U \delta q)$ corresponds to a positive Lyapunov exponent and then $(\delta q, \U \delta q)\notin G_-\cap G_+$. Then :  
  $$\begin{matrix}\hfill\frac{2 }{T}\log (\|\delta q(T)\|)+&\frac{\log 2C}{T}\geq \frac{1}{T}\log (\delta q (T)(\U-\stable ) (q_T, p_T)\delta q (T))\geq\hfill \\
 & \frac{1}{T}\log(\delta q (0)(\U-\stable)(q,p)\delta q (0))+\frac{1}{T}\int_0^Tm(\frac{\partial^2H}{\partial p^2}(q_t, p_t))q_+((\U-\stable)(q_t, p_t))dt.\end{matrix}$$
Using Birkhoff's ergodic theorem, we obtain: 
  $$\lambda(\mu)\geq\frac{1}{2}\int m(\frac{\partial^2H}{\partial p^2})q_+(\U-\stable)d\mu.$$
 \end{proof}
 \subsection{Twist maps: the weakly hyperbolic case}
\begin{theo}{\bf\ref{th4}}
Let $f: \A_d\rightarrow \A_d$ be a symplectic twist map and let $\mu$ be a locally minimizing ergodic measure with no zero Lyapunov exponents. We denote the smallest   Lyapunov exponent of $\mu$ by $\lambda (\mu)$ and  an upper bound for $\| s_1-s_{-1}\|$ above ${\rm supp} \mu$ by   $C$. Then we have: 
$$\lambda (\mu)\geq\frac{1}{2}\int \log\left( 1+\frac{1}{C}m(\U (x)-\esse (x))\right) d\mu (x).$$
\end{theo}

\demo We assume then that $(x_n)=(q_n, p_n)$ is a generic orbit for $\mu$.  Hence there exists $v_0\in G_+(x_0)\backslash \{ 0\}$ such that: 
$$\lim_{n\rightarrow \infty} \frac{1}{n}\log\left( \| Df^n(x_0)v_0\|\right)=\lambda (\mu).$$
The Lagrangian bundles $G_-$ and $G_+$ being transverse to the vertical at every point of ${\rm supp}\mu$, there exist two symmetric matrices $\esse$ and $\U$ such that $G_-$ (resp. $G_+$) is the graph of $\esse$ (resp. $\U$) in the usual coordinates of $\R^d\times\R^d=T_{x}\A_d$.
 As $G_-$ and $G_+$ are transverse $\mu$-almost everywhere,  we know that  there exists $\varepsilon>0$ such that $A_\varepsilon=\{ x\in {\rm supp}\mu; \U-\esse\geq \varepsilon {\bf 1}\}$ has positive $\mu$-measure. We may then assume that $x_0\in A_\varepsilon$ and that $\{ n; \U (x_n)-\esse(x_n)>\varepsilon{\bf 1}\}$ is infinite. Let us notice that in this case, $G_-$ and $G_+$ are transverse along the whole orbit of $x_0$ (but $\U-\esse$ can be very small at some points of this orbit).\\
Hence, for every $n\in\N$, there exists a unique positive definite matrix $S_0(x_n)$ such that: $S_0(x_n)^2=\U(x_n)-\esse (x_n)$. Let us recall that a matrix $M=\begin{pmatrix} a & b\\ c&d\\ \end{pmatrix}$ of dimension $2d$ is symplectic if and only if its entries satisfy the following equalities:
$${}^tac={}^tca;\quad {}^tbd={}^tdb;\quad {}^tda-{}^tbc={\bf 1}.$$
We define along the orbit of $x_0$ the following change of basis: 
$P=\begin{pmatrix} S_0^{-1}& S_0^{-1}\\ \esse S_0^{-1}& \U S_0^{-1}\\ \end{pmatrix}.$
Then it defines a symplectic change of coordinates, whose inverse is: 
$$Q=P^{-1}=\begin{pmatrix} 0&{\bf 1}\\ -{\bf 1}& 0\\ \end{pmatrix} {}^tP\begin{pmatrix} 0&-{\bf 1}\\ {\bf 1}& 0\\ \end{pmatrix}=\begin{pmatrix} S_0^{-1}\U & -S_0^{-1}\\ -S_0^{-1}\esse & S_0^{-1}\\ \end{pmatrix}.$$
We use this symplectic change of coordinates along the whole orbit of $x_0$. More precisely, if we denote the matrix of $Df^k$ in the usual canonical base $e=(e_i)$  by $M_k$,  then the matrix of $Df^k$ in the base $Pe=(Pe_i)$ is denoted by $\tilde M_k$; we have then: $\tilde M_k(x_n)= P^{-1}(x_{n+k})M_k(x_n)P(x_n)$. Let use notice that the image of the horizontal (resp. vertical) Lagrangian plane by $P$ is $G_-$ (resp. $G_+$). As the bundles $G_-$ and $G_+$ are invariant by $f$, we deduce that $\tilde M_k=\begin{pmatrix} \tilde a_k&0 \\
 0 &\tilde d_k\\ \end{pmatrix}$; we have  ${}^t\tilde a_k \tilde d_k={\bf 1}$ because this matrix is symplectic.\\
 Moreover, we know that: $M_k(x_n)=\begin{pmatrix} -b_k(x_n)s_{-k}(x_n)& b_k(x_n)\\
 c_k(x_n)& s_k(x_{k+n})b_k(x_n)\\ \end{pmatrix}$ where $G_k(x_n)=Df^k.V(x_{n-k})$ is the graph of $s_k(x_n)$.\\
 Writing that  $\tilde M_k(x_n)=\begin{pmatrix} \tilde a_k(x_n)&0 \\
 0 &\tilde d_k(x_n)\\ \end{pmatrix}= P^{-1}(x_{n+k})M_k(x_n)P(x_n)$, we obtain firstly: 
 $$S_0(x_{n+k})^{-1}{}^tb_k(x_n)S_0(x_n)^{-1}=S_0(x_{n+k})^{-1}(\esse (x_{n+k})-s_{k}(x_{n+k}))b_k(x_n)(s_{-k}(x_n)-\esse (x_n))S_0(x_n)^{-1};$$
 $$-S_0(x_{n+k})^{-1}{}^tb_k(x_n)S_0(x_n)^{-1}=S_0(x_{n+k})^{-1}(\U (x_{n+k})-s_{k}(x_{n+k}))b_k(x_n)(\U (x_n)-s_{-k}(x_n))S_0(x_n)^{-1}.$$
  We deduce that:
 $\tilde a_k(x_n)=S_0(x_{n+k})b_k(x_n)(\esse (x_n)-s_{-k}(x_n))S_0(x_n)^{-1}$  and: \\
  $\tilde d_k(x_n)=S_0(x_{n+k})b_k(x_n)(\U (x_n)-s_{-k}(x_n))S_0(x_n)^{-1}$.
  
 Because of the changes of basis that we used, $(\tilde a_k(x_n))_k$ represents the linearized dynamics $(Df^k_{|G_-(x_n)})_k$ restricted to $G_-$ and $(\tilde d_k(x_n))_k$ the linearized dynamics restricted to $G_+$.  Hence we need to study $(\tilde d_k(x_n))$ to obtain some information about the positive Lyapunov exponents of $\mu$. Let us compute:\\
 ${}^t\tilde d_k (x_n)=\tilde a_k(x_n)^{-1}=S_0(x_n)(\esse (x_n)-s_{-k}(x_n))^{-1}b_k(x_n)^{-1}S_0(x_{n+k})^{-1}$; we deduce:\\
 $${}^t\tilde d_k(x_n)\tilde d_k(x_n)=S_0(x_n)(\esse (x_n)-s_{-k}(x_n))^{-1}(\U (x_n)-s_{-k}(x_n))S_0(x_n)^{-1}$$
 $$=S_0(x_n)(\esse (x_n)-s_{-k}(x_n))^{-1}(\U (x_n)-\esse (x_n)+\esse (x_n)-s_{-k}(x_n))S_0(x_n)^{-1}$$
 $$={\bf 1}+ S_0(x_n)(\esse (x_n)-s_{-k}(x_n))^{-1} S_0(x_n)$$
 $$={\bf 1}+(\U (x_n)-\esse (x_n))^\frac{1}{2}(\esse (x_n)-s_{-k}(x_n))^{-1} (\U (x_n)-\esse (x_n))^\frac{1}{2}.$$
 Let us denote the conorm of $a$ (for the usual Euclidean norm of $\R^d$) by: $m(a)=\| a^{-1}\|^{-1}$. Then we have:
 $$m(\tilde d_k(x_n))^2=m({}^t\tilde d_k(x_n)\tilde d_k(x_n));$$
 and then: $m(\tilde d_k(x_n))^2\geq 1+\frac{1}{C}m((\U-\esse )(x_n))$ where $C$ designates $\sup \| s_1-s_{-1}\|$ above the (compact) support of $\mu$; indeed, we know that: $s_1-s_{-1}\geq \esse -s_{-k} >0$.
 
 The entry $\tilde d_k$ being multiplicative, we deduce that:
 $$m(\tilde d_k(x_0))^2\geq \prod_{n=0}^{k-1} (1+\frac{1}{C}m(\U (x_n)-\esse (x_n)))$$
 and: 
 $$\frac{1}{k}\log m(\tilde d_k(x_0))\geq \frac{1}{2k}\sum_{n=0}^{k-1}\log (1+\frac{1}{C}m(\U (x_n)-\esse (x_n))).$$
When $k$ tends to $+\infty$, we deduce from Birkhoff's ergodic theorem that: 
$$(*)\quad \liminf_{k\rightarrow \infty}\frac{1}{k}\log m(\tilde d_k(x_0))\geq\frac{1}{2}\int \log\left( 1+\frac{1}{C}m(\U (x)-\esse (x))\right) d\mu (x).$$
Let us recall that $(\tilde d_k(x_0))$ represents the dynamics along $G_+$, but the change of basis that we have done is not necessarily bounded. To obtain a true information about the Lyapunov positive exponents, we need to have a result for the matrix $D_k$ of $(Df^k_{|G_+(x_0)})$ in the base of $G_+$ whose matrix in the usual coordinates is: $\begin{pmatrix} {\bf 1}\\ \U \\ \end{pmatrix}$. Since $(\tilde d_k)$ is the matrix of $Df^k$ in the base whose matrix is $\begin{pmatrix} S_0^{-1}\\ \U S_0^{-1} \\ \end{pmatrix}$, we deduce that:  $D_k(x_0)=S_0(x_k)\tilde d_k(x_0)S_0(x_0)^{-1}$ and:\\
 $m(D_k(x_0))\geq m(S_0(x_k)) m(\tilde d_k(x_0)) m(S_0(x_0)^{-1})=\left( m(\U(x_k)-\esse (x_k))\right)^\frac{1}{2} m(\tilde d_k(x_0)) m(S(x_0)^{-1})$. \\
 We have $(*)$ and we know that: $\displaystyle{\liminf_{k\rightarrow \infty} m(\U (x_k)-\esse (x_k))\geq \varepsilon}$. We deduce: 
$$\lambda (\mu)\geq  \liminf_{k\rightarrow \infty}\frac{1}{k}\log m(D_k(x_0))\geq\frac{1}{2}\int \log\left( 1+\frac{1}{C}m(\U (x)-\esse (x))\right) d\mu (x).$$

\enddemo

\end{document}